\documentclass[
11pt,
article,
twoside,
openright,
a4paper]{memoir}

\usepackage[latin1]{inputenc}
\usepackage[T1]{fontenc}

\usepackage[english]{babel}

\usepackage{graphicx}

\usepackage[justification=RaggedRight]{subfig}

\usepackage{amsmath}
\usepackage{amssymb}
\usepackage{amsfonts}

\usepackage[amsmath,thmmarks,hyperref]{ntheorem}

\usepackage{mathtools}

\usepackage[sc]{mathpazo}

\makeatletter
\@ifpackageloaded{mathpazo}{}{%
  \usepackage{fix-cm}
  \usepackage{bbm}}
\makeatother

\usepackage[heavycircles]{stmaryrd}

\usepackage[ps,all,dvips]{xy}
\SelectTips{cm}{}
\CompileMatrices

\def\cxymatrix#1{\xy*[c]\xybox{\xymatrix#1}\endxy}

\usepackage{url}

\makeatletter
\@ifpackageloaded{mathpazo}{
  \newcommand{\mathsetfont}{\mathbb}
}{
  \newcommand{\mathsetfont}{\mathbbm}}
\makeatother

\newcommand{\DeclareMathSet}[1]{%
  \expandafter\newcommand\csname set#1\endcsname{\mathsetfont{#1}}}

\DeclareMathSet{C} 
\DeclareMathSet{R} 
\DeclareMathSet{Q} 
\DeclareMathSet{Z} 
\DeclareMathSet{N} 

\newcommand{\inv}{^{-1}}

\renewcommand{\phi}{\varphi}
\renewcommand{\epsilon}{\varepsilon}
\renewcommand{\rho}{\varrho}

\newcommand{\SL}{\mathrm{SL}} 

\newcommand{\ab}{\mathrm{ab}} 

\DeclareMathOperator{\Hom}{Hom}           
\DeclareMathOperator{\Aut}{Aut}           
\DeclareMathOperator{\Map}{Map}           
\DeclareMathOperator{\Fun}{Fun}           
\DeclareMathOperator{\Diff}{Diff}         
\newcommand{\boundary}{\partial}          
\DeclareMathOperator{\Ind}{Ind}           
\DeclareMathOperator{\Coind}{Coind}       

\makeatletter
\@ifpackageloaded{mathpazo}{

\renewcommand{\backslash}{\reflectbox{$/$}}
}{}
\makeatother

\newcommand{\sotimes}{\mathbin{\raise.14ex\hbox{$\scriptstyle\otimes$}}}

\ifpdf
\fi

\newcommand{\vchb}[1]{\vcenter{\hbox{#1}}}

\theoremnumbering{arabic} 
\theoremstyle{plain}      

\theoremsymbol{}

\theorembodyfont{\itshape} 
\theoremheaderfont{\normalfont\bfseries} 
\theoremseparator{.} 

\newtheorem{lemma}{Lemma}[chapter]
\newtheorem{theorem}[lemma]{Theorem}
\newtheorem{proposition}[lemma]{Proposition}
\newtheorem{corollary}[lemma]{Corollary}

\theorembodyfont{\upshape}

\newtheorem{remark}[lemma]{Remark}

\theoremstyle{nonumberplain}
\theoremheaderfont{\itshape}
\theorembodyfont{\normalfont}
\theoremsymbol{\ensuremath{\square}}
\newtheorem{proof}{Proof}

\theoremstyle{empty}

\qedsymbol{\ensuremath{\square}}

\setsecnumdepth{subsection}
\maxsecnumdepth{subsection}
\maxtocdepth{subsection}

\strictpagechecktrue

\title{Degree one cohomology with twisted coefficients of the mapping class group}
\author{Jørgen Ellegaard Andersen\and Rasmus Villemoes}

\begin{document}

\mainmatter

\maketitle

\begin{abstract}
  \noindent Let $\Gamma$ be the mapping class group of an oriented
  surface $\Sigma$ of genus $g$ with $r$ boundary components. We prove
  that the first cohomology group $H^1(\Gamma,
  \mathcal{O}(\mathcal{M}_{\SL_2(\setC)})^*)$ is non-trivial, where the
  coefficient module is the dual of the space of algebraic functions
  on the $\SL_2(\setC)$ moduli space over $\Sigma$.
\end{abstract}


\chapter{Introduction}
\label{cha:introduction}

Let $\Gamma = \Gamma_{g,r}$ denote the mapping class group of a
compact surface $\Sigma = \Sigma_{g,r}$ with genus $g$ and $r$
boundary components.  There is an action of $\Gamma$ on the moduli
space $\mathcal{M}_G$ of flat $G$-connections over $\Sigma$.
The vector space $\mathcal{O}(\mathcal{M}_G) \subseteq
\Fun(\mathcal{M}_G, \setC)$ of algebraic functions on the moduli space
is naturally a $\Gamma$-module. For the precise definition of the
class of algebraic functions we refer to the appendix.

Presently, we consider the special case of $G = \SL_2(\setC)$, and we
simply write $\mathcal{M}$ for $\mathcal{M}_{\SL_2(\setC)}$. In this
case, there is an isomorphism of $\Gamma$-modules
\begin{align}
  \label{eq:3}
  \nu\colon \mathcal{B}(\Sigma)\to\mathcal{O}(\mathcal{M}),
\end{align}
where the source denotes the algebra of BFK-diagrams on $\Sigma$: A
\emph{geometric BFK-diagram} on $\Sigma$ is a finite collection of
pairwise non-intersecting, non-trivial, unoriented simple loops on
$\Sigma$.  A BFK-diagram on $\Sigma$ is an isotopy class of geometric
BFK-diagrams.  Letting $B = B(\Sigma)$ denote the set of BFK-diagrams
on $\Sigma$, $\mathcal{B} = \mathcal{B}(\Sigma)$ is simply the complex
vector space spanned by $B$. There is a natural algebra structure on
this space; for details on this see~\cite{MR1691437}
and~\cite{ARS2006}. The isomorphism $\nu$ is given on a single simple
loop $\gamma$ by $\nu(\gamma) = -f_{\vec\gamma}$, where $\vec\gamma$
is any of the two oriented versions of $\gamma$, and $f_{\vec\gamma}$
is the function which to a gauge equivalence class $[A]$ of flat
connections associates the trace of the holonomy of $A$ along
$\vec\gamma$.


We may think of $\mathcal{B}$ as the set of maps $B\to\setC$ which
vanish except for a finite number of diagrams.  This is naturally
embedded in the larger module of \emph{all} maps $\hat{\mathcal{B}} =
\Map(B, \setC)$; this is clearly the same as the algebraic dual
$\mathcal{O}(\mathcal{M})^*$ of $\mathcal{O}(\mathcal{M})$. The action
of $\Gamma$ splits $B$ into orbits. Let $S$ denote a set of
representatives of these orbits, and for $D\in S$, let $\hat M_D$
(respectively $M_D$) denote the space of all maps from the orbit
through $D$ to $\setC$ (respectively, the maps $\Gamma D\to \setC$
which vanish for all but a finite number of diagrams in the orbit).
With this notation, we obtain splittings of $\mathcal{B}$ and
$\hat{\mathcal{B}}$ as $\Gamma$-modules
\begin{align}
  \mathcal{O}(\mathcal{M}) \cong \mathcal{B} &\cong \bigoplus_{D\in S}
  M_D \nonumber\\
  \label{eq:51}
  \mathcal{O}(\mathcal{M})^* \cong \hat{\mathcal{B}} &\cong
  \prod_{D\in S} \hat M_D
\end{align}
which induce decompositions in cohomology
\begin{align}
  \label{eq:52}
  H^*(\Gamma, \mathcal{B}) &\cong \bigoplus_{D\in S} H^*(\Gamma, M_D)\\
  \label{eq:53}
  H^*(\Gamma, \hat{\mathcal{B}}) &\cong \prod_{D\in S}
  H^*(\Gamma, \hat M_D).
\end{align}
A cocycle $u\colon \Gamma\to \mathcal{O}(\mathcal{M})^* =
\hat{\mathcal{B}} = \Map(B, \setC)$ may also be thought as a map
$u\colon \Gamma\times B\to\setC$ by simply putting $u(\gamma)(E) =
u(\gamma, E)$.
\begin{theorem}
  \label{thm:3}
  A cocycle $u\colon \Gamma\to \mathcal{O}(\mathcal{M})^* =
  \hat{\mathcal{B}} = \Map(B, \setC)$ is a coboundary if and only if
  for each $D\in S$, the restriction of $u$ to $\Gamma_D\times \{D\}$
  is identically $0$, where $\Gamma_D$ denotes the stabilizer of the
  diagram $D$ in $\Gamma$.
\end{theorem}
We will use this theorem to arrive at the main result:
\begin{theorem}
  \label{thm:5}
  For every $g, r\geq 0$, the cohomology group
  $H^1(\Gamma_{g,r}, \mathcal{O}(\mathcal{M})^*)$ is a direct product
  of summands $H^1(\Gamma, \hat M_D)$, each of which is
  finite-dimensional. Here $D$ runs over a set of representatives of
  BFK-diagrams on $\Sigma$.
\end{theorem}
In particular, we obtain by explicit examples
\begin{corollary}
  \label{cl:4}
  For $g\geq 1, r\geq 0$, $H^1(\Gamma_{g,r},
  \mathcal{O}(\mathcal{M})^*)$ is non-trivial.
\end{corollary}

The motivation to study the cohomology of the mapping class group with
these coefficients came from~\cite{0611126}, particularly
Proposition~6, where integrability of certain cocycles turn out to be
an obstruction to finding a $\Gamma$-invariant equivalence between two
equivalent star products on the moduli space. The motivation for
studying that problem comes from the expectation that the star
products discussed in~\cite{0611126} are equivalent to the star
product which is constructed in~\cite{MR1691437} and which is the same
as the ones induced on the $\SL_2(\setC)$-moduli space from the
constructions given in~\cite{MR1404925} and~\cite{MR1636568}.

This paper is organized as follows. In
Section~\ref{cha:group-cohom-backgr} we develop some of the basic
properties of group cohomology which are needed in the calculations,
ending with a proof of Theorem~\ref{thm:3}.  In
Section~\ref{cha:comp-h1gamma-mathc}, we develop an algorithm to
compute $H^1(\Gamma, \hat M_D)$ for any BFK-diagram $D$, which enables
us to prove Theorem~\ref{thm:5}. This is used in
Section~\ref{cha:an-example} to give a generic example of a
BFK-diagram for which the cohomology is non-zero. Finally we discuss
what we know when the coefficient module is
$\mathcal{O}(\mathcal{M})$.


\chapter{Group cohomological background}
\label{cha:group-cohom-backgr}

\begin{theorem}[Shapiro's Lemma]
  \label{thm:1}
  Let $H$ be a subgroup of $\Gamma$ and $A$ a left $H$-module. Then there
  are isomorphisms
  \begin{align}
    \label{eq:1}
    H_*(H, A) &\cong H_*(\Gamma, \Ind_H^\Gamma A)\\
    \label{eq:2}
    H^*(H, A) &\cong H^*(\Gamma, \Coind_H^\Gamma A).
  \end{align}
\end{theorem}
Here $\Ind_H^\Gamma$ is the so-called \emph{induced} module $\setZ
\Gamma \otimes_{\setZ H} A$, where $\setZ \Gamma$ is considered as a
right $H$-module via the right action of $H$ on $\Gamma$, and the left
$\Gamma$-module structure is given by $g\cdot(g' \sotimes a) = gg'
\sotimes a$ for $g,g'\in \Gamma$, $a\in A$. Similarly,
$\Coind_H^\Gamma A$ is the co-induced module $\Hom_{\setZ H}(\setZ
\Gamma, A)$ of $H$-equivariant maps from the left $H$-module $\setZ
\Gamma$ to $A$. The left action of $\Gamma$ is defined by
\begin{align*}
  (g\cdot f)(g') = f(g'g)
\end{align*}
for $g,g'\in \Gamma$, $f\in \Hom_{\setZ H}(\setZ \Gamma, A)$.

\begin{remark}
  \label{rem:1}
  If the action of $H$ on $A$ is trivial, there is a canonical
  bijection $\Hom_{\setZ H}(\setZ \Gamma, A) \to \Map(H\backslash
  \Gamma, A)$ given by $f\mapsto (Hg \mapsto f(g))$; equipping the
  latter with the $\Gamma$-action $(g\cdot f)(Hg') = f(Hg'g)$ this
  becomes an isomorphism of $\Gamma$-modules. The usual bijection
  between the sets of left and right cosets given by $Hg \mapsto g\inv
  H$ induces a bijection $\Map(H\backslash \Gamma, A)\to
  \Map(\Gamma/H, A)$, and the latter also carries a natural left
  $\Gamma$-action making this a $\Gamma$-isomorphism, namely $(g\cdot
  f)(g'H) = f(g\inv g'H)$.
\end{remark}
We summarize the special case of Shapiro's Lemma we will need in a
corollary:
\begin{corollary}
  \label{cl:3}
  Let $A$ be an abelian group, and $\Gamma$ a group which acts
  transitively on a set $R$. Consider the $\Gamma$-module $\Map(R, A)$
  of all maps $R\to A$ with action given by $(g\cdot f)(r) = f(g\inv
  r)$. Let $D\in R$ be any element, and $\Gamma_{D}\subseteq \Gamma$
  the stabilizer subgroup of $D$. Then there is an isomorphism
  \begin{align}
    \label{eq:63}
    H^*(\Gamma, \Map(R, A)) \cong H^*(\Gamma_D, A)
  \end{align}
  where $A$ is considered as a trivial $\Gamma_D$-module.
\end{corollary}
\begin{proof}
  The bijection $\Gamma/\Gamma_D\to R$ given by $g\Gamma_D\mapsto
  gD$ clearly induces an isomorphism of $\Gamma$-modules
  $\Map(\Gamma/\Gamma_D, A)\to \Map(R, A)$. Then from Shapiro's Lemma
  and the isomorphisms mentioned in the above remark we have a
  sequence of isomorphisms
  \begin{align*}
    H^*(\Gamma_D, A) &\cong H^*(\Gamma, \Hom_{\setZ \Gamma_D}(\setZ \Gamma, A))\\
    &\cong H^*(\Gamma, \Map(\Gamma_D\backslash \Gamma, A))\\
    &\cong H^*(\Gamma, \Map(\Gamma/\Gamma_D, A))\\
    &\cong H^*(\Gamma, \Map(R, A)).
  \end{align*}
\end{proof}

Note that the $\Gamma$-module $\Map(R, A)$ can also be considered as the
set of all formal $A$-linear combinations of elements from $R$ (that
is, the sum $\sum_{r\in R} m_r r$ corresponds to the map $r\mapsto
m_r$).

Specializing to the case $*=1$, we will now describe a more or less
explicit isomorphism $H^1(\Gamma, \Map(R, A))\to H^1(\Gamma_D, A)$.
First note that a map $u\colon \Gamma\to \Map(R, A)$ can equally well
be considered as a map $u\colon \Gamma\times R\to A$ by the adjoint
formula $u(g)(r) = u(g, r)$. In this context, the cocycle condition
reads
\begin{align}
  \label{eq:28}
    u(g_1g_2, r) = u(g_1, r) + u(g_2, g_1\inv r).
\end{align}
We wish to derive necessary and sufficient conditions for a cocycle
$u$ to be a coboundary $\delta f$. For the rest of this section, fix
an element $D\in R$ and let $\Gamma_{D}\subseteq \Gamma$ denote the
stabilizer subgroup of $D$.
\begin{lemma}
  \label{lem:1}
  A cocycle $u\colon \Gamma\times R\to A$ is a coboundary if and only
  if, for every pair $g_1, g_2\in \Gamma$ with $g_1 g_2\inv \in \Gamma_D$,
  $u$ satisfies the condition
  \begin{align}
    \label{eq:33}
    u(g_1, D) = u(g_2, D).
  \end{align}
\end{lemma}
\begin{proof}
  First we prove the necessity of the condition. Suppose that $u =
  \delta f$ for some $f\colon R\to A$. Since the action is transitive,
  it is easy to see that the kernel of $\delta\colon C^0(\Gamma,
  \Map(R, A)) \to C^1(\Gamma, \Map(R, A))$ is the set of constant maps
  $R\to A$. Thus we may WLOG assume that $f(D) = 0$. Recall that $u =
  \delta f$ means that for every $g\in \Gamma$, $r\in R$ we have $u(g,
  r) = f(r) - f(g\inv r)$. In particular,
  \begin{align}
    \label{eq:29}
    f(g\inv D) = -u(g, D)
  \end{align}
  Now if $g_1g_2\inv \in \Gamma_{D}$, we have $g_1\inv D = g_2\inv
  D$, and thus $-u(g_1, D) = f(g_1\inv D) = f(g_2\inv D) =
  -u(g_2, D)$ as desired.

  Now suppose that $u$ satisfies~\eqref{eq:33} whenever $g_1g_2\inv
  D = D$. We need to construct a map $f\colon R\to A$. For
  $r\in R$, choose $g\in \Gamma$ so that $g\inv D = r$, and define $f$
  using \eqref{eq:29}, ie. $f(r) = f(g\inv D) = -u(g, D)$. By
  assumption, this is a well-defined map (independent of the chosen
  $g$), and we only need to check that $u = \delta f$. Let $h\in \Gamma$
  and $r\in R$ be arbitrary. To calculate $(\delta f)(h, r)$, we may
  choose any $g\in \Gamma$ with $g\inv D = r$, and we obtain
  \begin{align*}
    (\delta f)(h, r) = f(r) - f(h\inv r) &= f(g\inv D) - f( (gh)\inv
    D)\\
    &= -u(g, D) + u(gh, D) = u(h, g\inv D) = u(h, r)
  \end{align*}
  by the cocycle condition \eqref{eq:28}.
\end{proof}
\begin{lemma}
  \label{lem:2}
  The restriction of $u$ to $\Gamma_D\times \{D\}$ is a group
  homomorphism $\tilde u\colon \Gamma_D \to A$.
\end{lemma}
\begin{proof}
  Let $g, h\in \Gamma_D$. Then
  \begin{align}
    \tilde u(gh) = u(gh, D) &= u(g, D) + u(h, g\inv
    D) \nonumber\\
    \label{eq:30}
    &= u(g, D) + u(h, D) = \tilde u(g) + \tilde u(h)
  \end{align}
  as claimed.
\end{proof}
Since $A$ is abelian, $\tilde u$ factors through the abelinization
$(\Gamma_D)_{\ab}$ of $\Gamma_D$, and we have thus established a map
$\phi\colon Z^1(\Gamma, \Map(R, A))\to \Hom(\Gamma_D, A) =
\Hom((\Gamma_D)_{\ab}, A)$.  The latter group may be thought of as the
cohomology group $H^1((\Gamma_D)_{\ab}, A)$ with trivial action of
$(\Gamma_D)_{\ab}$ on $A$.

\begin{theorem}
  \label{thm:2}
  The map $\phi$ factors to an isomorphism $H^1(\Gamma, \Map(R, A))\to
  H^1((\Gamma_D)_{\ab}, A)$, which is also denoted $\phi$.
\end{theorem}
Before we begin the proof, we need an observation: For any cocycle $u$
and any $g\in \Gamma$, $h\in \Gamma_D$ we have
\begin{align*}
  u(ghg\inv,  gD) &= u(g, gD) + u(hg\inv, D)\\
  &= u(g, gD) + u(h, D) + u(g\inv, D)\\
  &= u(h, D)
\end{align*}
using $h\inv D = D$ and the fact that $0 = u(1) = u(g\inv\cdot g) =
u(g\inv) + g\inv.u(g)$.
\begin{proof}[of Theorem~\ref{thm:2}]
  To prove the first part of the theorem, we need to show that the
  restriction of a cobundary $\delta f$ to $\Gamma_D\times \{D\}$ is
  identically $0$. But this is trivial since
  \begin{align*}
    \widetilde{\delta f}(h) = (\delta f)(h, D) = f(D) - f(h\inv D) = 0
  \end{align*}
  for $h\in \Gamma_D$.

  Next, assume that the cocycle $u$ restricts to the zero homomorphism
  $\Gamma_D\to A$. Then for any two elements $g_1, g_2\in \Gamma$ with
  $g_1g_2\inv \in \Gamma_D$ we have
  \begin{align*}
    0 &= u(g_1 g_2\inv, D)\\
    &= u(g_1, D) + u(g_2\inv, g_1\inv D)\\
    &= u(g_1, D) + u(g_2\inv)(g_1\inv D)\\
    &= u(g_1, D) - g_2\inv.u(g_2)(g_1\inv D)\\
    &= u(g_1, D) - u(g_2)(g_2 g_1\inv D)\\
    &= u(g_1, D) - u(g_2, D)
  \end{align*}
  since $g_2g_1\inv = (g_1g_2\inv)\inv \in \Gamma_D$, and by
  Lemma~\ref{lem:1} we see that $u$ is a coboundary. This shows that
  $\phi$ is injective.

  Now, for surjectivity, let $u\colon \Gamma_D\to A$ be any
  homomorphism. We need to extend $u$ to all of $\Gamma\times R$ in
  such a way that it becomes a cocycle. To produce this extension, we
  first assume that an extension exists, and use this to write a
  formula for a cocycle cohomologous to the given extension. Then we
  prove that this formula actually defines a cocycle.

  Choose a collection $\{h_i\}_{i\in I}$ of representatives for the
  set $\Gamma_D\backslash \Gamma$ of \emph{right} cosets of
  $\Gamma_D$, and let $1\in \Gamma$ represent the coset $\Gamma_D$.
  Recall that the map $\Gamma_D\backslash \Gamma \to \Gamma/\Gamma_D$
  given by $\Gamma_Dx \mapsto x\inv \Gamma_D$ is a bijection between
  the set of right cosets and the set of left cosets of $\Gamma_D$. In
  particular, $\{h_i\inv\}_{i\in I}$ is a collection of
  representatives of the set of left cosets. We also have a bijection
  $\Gamma/\Gamma_D\to R$ given by $x\Gamma_D\mapsto xD$. Now, for any
  coboundary $\delta f$ with which we alter $u$, we may (as has been
  used a couple of times) WLOG assume that $f(D) = 0$. Then the
  formula $(\delta f)(h_i)(D) = f(D) - f(h_i\inv D) = -f(h_i\inv D)$
  and the fact that $i\mapsto h_i\inv D$ is a bijection $I\to R$ show
  that we may assume that the extension $u$ satisfies $u(h_i, D) = 0$
  for $i\in I$.  Furthermore, $u$ is uniquely determined by its
  cohomology class and this requirement.

  The cocycle condition implies that
  \begin{align}
    \label{eq:36}
    u(gh_i, D) = u(g, D) + u(h_i, g\inv D) = u(g, D)
  \end{align}
  for $i\in I$ and $g\in \Gamma_D$. Since every $x\in \Gamma$ admits a
  unique factorization as $x = gh_i$ for some $i\in I$ and $g\in
  \Gamma_D$, this formula extends $u$ to all of $\Gamma\times \{D\}$.

  Now consider any $x\in \Gamma$ and $E\in R$. There is a unique $j\in
  I$ with $h_j\inv D = E$, and we have $\Gamma_E = h_j\inv \Gamma_D
  h_j$. Furthermore, the collection $\{h_j\inv h_i h_j\}_{i\in I}$ is
  a collection of representatives for the set $\Gamma_E\backslash
  \Gamma$ of right cosets of $\Gamma_E$. This means that we may
  factorize $x$ uniquely as $(h_j\inv g_0 h_j)(h_j\inv h_i h_j)$ for
  some $g_0\in \Gamma_D$ and $i\in I$. Now we calculate
  \begin{align}
    \label{eq:37}
    u(x, E) &= u(h_j\inv g_0 h_j\cdot h_j\inv h_i h_j, h_j\inv D)\\
    &= u(h_j\inv g_0 h_j, h_j\inv D) + u(h_j\inv h_i h_j, h_j\inv
    g_0\inv h_j h_j\inv D)
  \end{align}
  By the observation preceding this proof (with $g=h_j\inv$ and
  $h=g_0$), the first term is equal to the known quantity $u(g_0,
  D)$. For the second term, we apply the cocycle condition a few
  more times:
  \begin{align*}
    u(h_j\inv h_i h_j, h_j\inv D) &= u(h_j\inv, h_j\inv D) + u(h_i
    h_j, D)\\
    &= -u(h_j, D) + u(h_i h_j, D)\\
    &= u(h_i h_j, D)
  \end{align*}
  which is also known since $u$ is known on $\Gamma\times \{D\}$%
  . Thus our formula for the
  extension of $u$ to all of $\Gamma\times R$ reads
  \begin{align}
    \label{eq:38}
    u(x, E) = u(g_0, D) + u(h_i h_j, D)
  \end{align}
  where $j\in I$ is the unique index such that $h_j\inv D = E$,
  $i\in I$ is the unique index so that $x$ belongs to the right coset
  of $\Gamma_E$ represented by $h_j\inv h_i h_j$, and $g_0 = h_j g h_j\inv$
  is the unique element i $\Gamma_D$ such that $x = g(h_j\inv h_i h_j) =
  (h_j\inv g_0 h_j)(h_j\inv h_i h_j)$. The second term above is
  defined by \eqref{eq:36}; thus one must find the $k\in I$ such that
  $h_i h_j$ is an element of the right coset of $\Gamma_D$ represented by
  $h_k$, say $h_ih_j = g_1h_k$ for $g_1\in \Gamma_D$, and then $u(h_i h_j,
  D) = u(g_1, D)$. It remains to check that \eqref{eq:38} defines
  a cocycle.

  Let $x,y\in \Gamma$ and $E\in R$ be arbitrary. As above, there is a
  unique $j\in I$ with $h_j\inv D = E$. Lets try to calculate the
  right-hand side of the cocycle condition $u(xy, E) = u(x, E) + u(y,
  x\inv E)$. We must choose $i\in I$ and $g_1\in \Gamma_D$ such that
  \begin{align}
    \label{eq:39}
    x &= (h_j\inv g_1 h_j)(h_j\inv h_i h_j)
  \end{align}
  and next we choose $k\in I$ and $g_2\in \Gamma_D$ such that $h_i h_j
  = g_2 h_k$. Then
  \begin{align*}
    u(x, E) = u(g_1, D) + u(g_2, D) = u(g_1g_2, D)
  \end{align*}
  Now, the element $x\inv E$ of $R$ is the same as
  \begin{align*}
    x\inv E = h_j\inv h_i\inv g_1\inv h_j E = h_j\inv h_i\inv D =
    (h_i h_j)\inv D = (g_2 h_k)\inv D = h_k\inv D
  \end{align*}
  so in the calculation of $u(y, x\inv E)$ it is $h_k$ which plays the
  role as $h_j$ in the recipe. This recipe then requires us to find
  $g_3\in \Gamma_D$ and $\ell \in I$ such that
  \begin{align}
    \label{eq:40}
    y = (h_k\inv g_3 h_k)(h_k\inv h_\ell h_k),
  \end{align}
  and $g_4\in \Gamma_D$ and $m\in I$ such that $h_\ell h_k = g_4 h_m$. Then
  \begin{align*}
    u(y, x\inv E) = u(g_3, D) + u(g_4, D) = u(g_3g_4, D).
  \end{align*}

  Multiplying $x$ and $y$ using the expressions \eqref{eq:39} and
  \eqref{eq:40} and the relations defining the various $h$'es we
  obtain
  \begin{align}
    \label{eq:45}
    xy &= (h_j\inv g_1 h_i h_j)(h_k\inv g_3 h_\ell h_k) \nonumber\\
    &= h_j\inv g_1 g_2 g_3 g_4 h_m
  \end{align}

  On the other hand, the recipe requires us to choose $g\in \Gamma_D$
  and $n\in I$ such that
  \begin{align}
    \label{eq:42}
    xy &= h_j\inv g h_j h_j\inv h_n h_j,
  \end{align}
  and $g'\in \Gamma_D$ and $p\in I$ such that $h_n h_j = g' h_p$. Then
  $u(xy, E) = u(g, D) + u(g', D)$. Comparing \eqref{eq:45} and
  \eqref{eq:42} we see that $g_1 g_2 g_3 g_4 h_m = g h_n h_j$, showing
  that (by uniqueness of $g'$ and $p$) $h_p = h_m$ and
  \begin{align}
    \label{eq:43}
    g' &= g\inv g_1g_2g_3g_4
  \end{align}
  Finally we conclude that
  \begin{align*}
    u(xy, E) &= u(g, D) + u(g', D)\\
    &= u(g_1g_2g_3g_4, D)\\
    &= u(g_1g_2, D) + u(g_3g_4, D)\\
    &= u(x, E) + u(y, x\inv E)
  \end{align*}
  showing that the given recipe in fact defines a cocycle $u\colon
  \Gamma\times R\to A$. The proof is complete.
\end{proof}
\begin{proof}[of Theorem~\ref{thm:3}]
  By the splitting~\eqref{eq:53}, a cocycle $u\colon\Gamma\to
  \hat{\mathcal{B}}$ is the same as a collection of cocycles
  $u_D\colon \Gamma\to \hat M_D$ for $D\in S$. In fact, thinking of
  $u$ as a map $\Gamma\times B\to \setC$, $u_D$ is simply the
  restriction of $u$ to $\Gamma\times (\Gamma D)$. Specializing
  Theorem~\ref{thm:2} to the case $A = \setC$ and $R = \Gamma D$, we
  see that each $u_D$ is a coboundary if and only if $u_D$ restricted
  to $\Gamma_D \times \{D\}$ is zero.
\end{proof}

In section~\ref{cha:comp-h1gamma-mathc} below, we are going to need a
theorem linking the low-dimensional homology groups of the groups
appearing in a short exact sequence. Again quoting from
\cite{MR672956} (Corollary~VII.6.4)
\begin{theorem}
  \label{thm:4}
  Let $1\to A\to B\to C\to 1$ be a short exact sequence of groups, and $M$ a
  $B$-module. Then there is an exact sequence of low-dimensional homology groups
  \begin{align}
    \label{eq:60}
    \begin{split}
      H_2(B, M) \to H_2(C, M_A) &\to\\
      H_1(A, M)_C \to H_1(B, M) \to H_1(C, M_A) &\to 0.
    \end{split}
  \end{align}
\end{theorem}
Here we regard $M$ as an $A$-module via restriction of scalars, and
then clearly $C \cong B/A$ acts on the co-invariants group $M_A$,
making sense of $H_*(C, M_A)$. Since $A$ is a normal in $B$,
conjugation by $b\in B$ defines an action on $A$ by automorphisms, so
there is an induced action on homology $c(b)_*\colon H_*(A, M)\to
H_*(A, M)$. One may show that $A$ acts trivially, so there is an
induced action of $C$, and we have $H_1(A, M)_B = H_1(A, M)_C$.

\chapter{Computing $H^1(\Gamma, \mathcal{O}(\mathcal{M})^*)$}
\label{cha:comp-h1gamma-mathc}

By the $\Gamma$-equivariant isomorphism~\eqref{eq:3} and the
splitting~\eqref{eq:53}, it is clear that $H^1(\Gamma,
\mathcal{O}(\mathcal{M})^*)$ splits as a direct product of
$H^1(\Gamma, \hat M_D)$-s, proving the first part of
Theorem~\ref{thm:5}. In order to prove that these are all
finite-dimensional, we develop in this section an algorithm to compute
them. Then it suffices to find a single diagram $D$ for which
$H^1(\Gamma, \hat M_D)$ is non-zero in order to prove
Corollary~\ref{cl:4}.

Recall that $\hat M_D = \Map(\Gamma D, \setC)$.  By the previous
section (specifically Corollary~\ref{cl:3}), we have $H^1(\Gamma, \hat
M_D) \cong H^1(\Gamma_D, \setC) = \Hom(\Gamma_D, \setC)$, so since
$\setC$ is abelian and torsion-free, to compute $H^1(\Gamma, \hat
M_D)$ amounts to computing the first homology group of the stabilizer
$\Gamma_D$ with rational coefficients.

In order to compute $H_1(\Gamma_D, \setQ)$, we consider the surface
$\Sigma'$ which is obtained from $\Sigma$ by cutting along $D$.
Let $n$ denote the number of components of $D$, and let $n'$ be the
maximal number of components of $D$ such that $\Sigma$ cut along these
is still connected. Put $n = n'+n''$. Then $\Sigma'$ is a (possibly
non-connected) surface with $1+n''$ connected components, total genus
$g' = g - n'$ and a total of $r' = r + 2n$ boundary components. There
is a ``glueing map'' $j\colon \Sigma'\to\Sigma$ which is a local
diffeomorphism away from the $2n$ boundary components arising from
$D$. The mapping class group $\Gamma'$ of $\Sigma'$ maps to $\Gamma_D$
(via $j$), because for any homeomorphism $\gamma'\colon
\Sigma'\to\Sigma'$ fixed on the boundary $\boundary \Sigma'$, there is
a unique homeomorphism $\gamma\colon\Sigma\to\Sigma$ fixing
$\boundary\Sigma$ and fitting into the diagram
\begin{align}
  \label{eq:58}
  \cxymatrix{{ \Sigma' \ar[r]^{\gamma'} \ar[d]^j & \Sigma' \ar[d]^j\\
    \Sigma \ar[r]^{\gamma} & \Sigma }}
\end{align}
and clearly any isotopy fixed on $\boundary \Sigma'$ descends to an
isotopy fixed on $\boundary \Sigma$.
This group homomorphism is not injective, since a diffeomorphism of
$\Sigma'$ consisting of ``oppositely oriented'' Dehn twists along two
boundary components glued together by $j$ is isotopic to the identity
in $\Diff(\Sigma)$. Also, it is not surjective, since the elements of
$\Gamma_D$ are allowed to permute the components of~$D$, which no
homeomorphism coming from $\Gamma'$ can do.

Hence, we need a notion of a ``larger'' mapping class group through
which information about $H_*(\Gamma')$ can be translated into
information about $H_*(\Gamma_D)$. To this end, choose an oriented
parametrization $c\colon \bigsqcup_{2n+r} S^1 \to \boundary \Sigma'$
of the boundary.  Consider the group $\Diff(\Sigma'; c)$ of
diffeomorphisms of $\Sigma'$ preserving this parametrization, ie. the
group of diffeomorphisms $\gamma$ such that $c\inv \circ
\smash{\gamma_{\vert\boundary\Sigma'}}\circ c$ is a permutation of the
$2n+r$ copies of $S^1$ consisting of identity maps. We define $G'$ to
be the group $\Diff(\Sigma'; c)$ modulo isotopies preserving $c$. Note
that in case $\Sigma'$ is non-connected, elements of $G'$ are allowed
to permute homeomorphic components.

There is a homomorphism from $G'$ to the permutation group of the set
of boundary components of $\Sigma'$, and the kernel of this map is
easily seen to be $\Gamma'$ (if a diffeomorphism maps each boundary
component to itself and at the same time preserves a parametrization,
it fixes the boundary point-wise). Thus we have a short exact sequence
\begin{align*}
  \xymatrix{ 1\ar[r] & \Gamma' \ar[r] & G' \ar[r] & P' \ar[r] & 1,}
\end{align*}
where $P'$ is the appropriate subgroup of $S_{\pi_0 \boundary
  \Sigma'}$. (In case $\Sigma'$ is connected, any permutation of the
boundary components is realizable through a diffeomorphism.)

Consider the subgroup $Q'\subseteq P'$ of permutations which fix
$\pi_0(\boundary \Sigma)$ (or more precisely the set
$\pi_0(j\inv(\boundary \Sigma))$) and preserves the pairing of
elements of $\pi_0(\boundary \Sigma' - j\inv(\boundary \Sigma))$
induced by $j$. In other words, $Q'$ consists of the permutations
$\sigma\in P'$ such that $\sigma(b) \in \pi_0(j\inv(j(b))$ for every
boundary component~$b$ of~$\Sigma'$; ie.~$\sigma(b) = b$ if $b$ is a
boundary component of $\Sigma$, otherwise $\sigma(b)$ is either equal
to $b$ or the boundary component of $\Sigma'$ which it is identified
with by $j$. Let $H'$ be the pre-image of $Q'$ in $G'$, so we have a
new exact sequence
\begin{align}
  \label{eq:57}
  \xymatrix{ 1\ar[r] & \Gamma' \ar[r] & H' \ar[r] & Q' \ar[r] & 1.}
\end{align}
One could also define $H'$ as the subgroup of $G'$ consisting of
elements which descend to elements of $\Gamma_D$ as in \eqref{eq:58}
above. Now, the homomorphism $H'\to \Gamma_D$ is easily seen to be
surjective. For the moment, assume that the kernel of this map is the
free abelian group $\setZ^n$ with one generator for each component of
$D$. Then we have another short exact sequence
\begin{align}
  \label{eq:59}
  \xymatrix{ 1 \ar[r] & {\setZ^n} \ar[r] & H' \ar[r] & \Gamma_D \ar[r]
    & 1.}
\end{align}

Lets apply Theorem~\ref{thm:4} to \eqref{eq:57}. As explained earlier,
we use rational coefficients (with trivial action), so since $Q'$ is a
finite group, its rational homology (in positive dimensions) vanishes,
and we are left with
\begin{align}
  \label{eq:61}
  \xymatrix{ 0 \ar[r] & H_1(\Gamma'; \setQ)_{Q'} \ar[r] & H_1(H';
    \setQ) \ar[r] & 0}
\end{align}

By work of (among others) Harer, the low-dimensional homology groups
of mapping class groups are known, at least with rational
coefficients. For easy reference, we collect the results we will need
in a proposition.
\begin{proposition}
  \label{prop:2}
  Let $\Gamma_{g,r}$ denote the mapping class group of a genus $g$
  surface with $r$ boundary components.
  \begin{enumerate}[\normalfont(1)]
    \firmlist
  \item If $g\geq 2$, $H_1(\Gamma_{g,r}; \setQ) = 0$.
    \label{genus-geq-2}
  \item For $g = 1$, $H_1(\Gamma_{1, r}; \setQ) \cong \setQ^{r}$. \label{genus-eq-1}
  \item For $g = 0$, $H_1(\Gamma_{0, r}; \setQ) \cong
    \setQ^{(r-1)r/2}$. \label{genus-eq-0}
  \end{enumerate}
\end{proposition}
In fact,~(\ref{genus-geq-2}) is easy to prove using simple geometric
considerations, and the fact that the mapping class group is generated
by Dehn twists (cf. the appendix, Corollary~\ref{cl:1}), except that
in genus $2$ one has to rely on a presentation of the mapping class
group. For proofs of (\ref{genus-eq-1}) and (\ref{genus-eq-0}) we
refer to~\cite{MR1106936}. Harer's work even give explicit generators:
In the case~(\ref{genus-eq-1}), the Dehn twists along the $r$ boundary
components represent a basis for the rational homology. In the
case~(\ref{genus-eq-0}), think of $\Sigma_{0,r}$ as the closed unit
disc with $r-1$ small open discs centered at the $x$-axis removed.
Then the $r-1$ Dehn twists along the boundaries of these small discs,
along with $\binom{r-1}{2} = (r-2)(r-1)/2$ twists along circles
enclosing exactly two of these discs represent a basis for
$H_1(\Gamma_{0, r}; \setQ)$.

Since the mapping class group of a non-connected surface, where each
connected component has at least one boundary component, is obviously
the product of the mapping class groups of the components, this shows
how to find a set of generators for $H_1(\Gamma'; \setQ)$, and that we
may in fact represent these generators by Dehn twists. Since the
action of $Q'$ is induced by the conjugation action of $H'$ on
$\Gamma'$, we see 
that the action simply
identifies some of these generators. The exact details regarding which
generators are thus identified depend on topological constraints (for
instance, it can happen that two components of $\Sigma'$ are
homeomorphic, but that there does not exist a glueing-compatible
diffeomorphism taking one to the other).

Applying Theorem~\ref{thm:4} to \eqref{eq:59} and using the surjective
map $H_1(A, M)\to H_1(A, M)_C$, we obtain another exact sequence
\begin{align*}
  \xymatrix{ H_1(\setZ^n; \setQ) \ar[r] & H_1(H'; \setQ) \ar[r] &
    H_1(\Gamma_D; \setQ) \ar[r] & 0, }
\end{align*}
which by using the isomorphism~\eqref{eq:61} becomes
\begin{align}
  \label{eq:31}
  \xymatrix{ H_1(\setZ^n; \setQ) \ar[r] & H_1(\Gamma'; \setQ)_{Q'} \ar[r] &
    H_1(\Gamma_D; \setQ) \ar[r] & 0. }
\end{align}
Recall that the group $\setZ^n$ really means the free abelian group
generated by~$n$ pairs of left and right Dehn twists in the boundary
components arising from the cutting along $D$. Consider such a pair of
Dehn twists $\tau_b$, $\tau_{b'}\inv$ in the boundary curves
$b, b'$ which are glued together by $j$. If both $b$ and $b'$ belong
to components of $\Sigma'$ with genus $\geq 2$, $\tau_b\tau_{b'}\inv$
is mapped to $0$ in $H_1(\Gamma'; \setQ)_{Q'}$. If exactly one of them
belongs to a component with genus $\geq 2$, the product
$\tau_b\tau_{b'}\inv$ is mapped to a generator of $H_1(\Gamma';
\setQ)_{Q'}$, and if both $b,b'$ belong to genus $\leq 1$ components,
the image of $\tau_b\tau_{b'}\inv$ in $H_1(\Gamma'; \setQ)_{Q'}$
identifies the generators $\tau_b$,~$\tau_{b'}$ (in case these were
not already identified by the action of~$Q'$).

In this way we see that we have a combinatorial method to compute
$H_1(H; \setQ)$ by cutting along $D$, writing down all generators coming
from genus~$0$ and~$1$ components, and identifying and/or removing
generators according to which permutations of the generators are
topologically realizable or which are killed by the image of
$H_1(\setZ^n;\setQ)$.
\begin{proposition}
  \label{prop:1}
  For a BFK-diagram $D$ on $\Sigma$, $H^1(\Gamma, \hat M_D) \cong
  \setC^n$, where $n$ is the dimension of the rational vector space
  $H_1(\Gamma_D; \setQ)$. It is always finite, and may be found by the
  algorithm described above.
\end{proposition}
This in particular proves the second claim in Theorem~\ref{thm:5}.

\chapter{An example}
\label{cha:an-example}

We now wish to construct a diagram $D$ such that the factor
$H^1(\Gamma, \hat M_D)$ in \eqref{eq:53} is non-trivial, ie. such that
there exists a non-zero homomorphism $\Gamma_D\to \setC$. In view
of~\eqref{eq:31} and Proposition~\ref{prop:2} above, we must choose
$D$ such that $\Sigma'$ contains at least one component of genus at
most $1$.
With this in mind, we arrive at the following generic example.

\begin{figure}[hbt]
  \centering
  \includegraphics{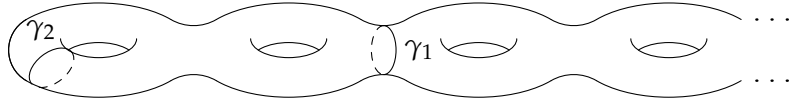}
  \caption{A two-component BFK-diagram.}
  \label{fig:1}
\end{figure}
Consider a surface $\Sigma_{g,r}$ of genus $g\geq 4$ with an arbitrary
number, $r$, of boundary components. We may choose a two-component
BFK-diagram $D = \gamma_1 \cup \gamma_2$ such that $\gamma_1$ is a
seperating curve, dividing $\Sigma$ into a surface of genus $2$ with
one boundary component (arising from the cut), and a surface of genus
$g-2$ containing all the original boundary components. The curve
$\gamma_2$ is chosen such that it is a non-seperating (hence
genus-decreasing) curve in the genus $2$ component. Hence in this
case, the cut surface $\Sigma'$ consists of a genus $1$ component
$\Sigma'_1$ with three boundary components (one arising from the cut
along $\gamma_1$, the other two being the boundaries arising from the
cut along $\gamma_2$) and a genus $g-2$ component $\Sigma'_2$ with
$r+1$ boundary components. We denote the boundary curves arising from
the cut along $\gamma_1$ by $\eta_{11}$ and $\eta_{12}$, respectively
(with $\eta_{11}$ belonging to $\Sigma'_1$), and the boundaries
arising from the cut along $\gamma_2$ are similarly denoted
$\eta_{21}$, $\eta_{22}$.

Now clearly the mapping class group $\Gamma'$ is the product
$\Gamma_{g-2, r+1}\times \Gamma_{1, 3}$, and since $g\geq 4$, the
rational homology of the first factor vanishes. By
Proposition~\ref{prop:2}, $H_1(\Gamma'; \setQ) = H_1(\Gamma_{1,3};
\setQ) = \setQ\{\eta_{11}, \eta_{21}, \eta_{22}\}$, the rational
vector space spanned by the Dehn twists $\eta_{11}, \eta_{21},
\eta_{22}$ (in order to keep notation simple we denote a curve and the
Dehn twist along it by the same symbol). It is clear that the only
(non-trivial) permutation of the boundary components of $\Sigma'$ that
is compatible with glueing is the interchange of $\eta_{21}$ and
$\eta_{22}$, and the effect of $Q' \cong \setZ/2$ on $H_1(\Gamma';
\setQ)$ is hence to identify the two generators $\eta_{21},
\eta_{22}$, so we may think of $H_1(\Gamma'; \setQ)_{Q'}$ as the
rational vector space spanned by $\eta_{11}$ and
$\eta_{21}$.

The two generators of $H_1(\setZ^2; \setQ)$ may be chosen to be
$\eta_{11} - \eta_{12}$ and $\eta_{21} - \eta_{22}$. We have just
established that $\eta_{21} = \eta_{22}$ in $H_1(\Gamma';
\setQ)_{Q'}$, so the latter of these generators is mapped to $0$. We
also have that $\eta_{12} = 0$ in $H_1(\Gamma'; \setQ)_{Q'}$, so the
first generator is simply mapped to $\eta_{11}$. Hence by the
exactness of \eqref{eq:31}, $H_1(\Gamma_D; \setQ)$ is a
$1$-dimensional rational vector space spanned by the Dehn twist in
$\gamma_2$, so for this particular diagram $D$, we see that
$H^1(\Gamma, \hat M_D)\cong \setC$.

It is not hard to obtain similar examples for the remaining low values
of genus. If $g = 3$, we may choose $D = \gamma_1\cup \gamma_2$ as
above, dividing the surface into two genus $1$ components, one with three
boundary components arising from the cut along $D$ and one with the
original boundary components (if any) together with the one arising
from the cut along $\gamma_1$. The only difference is that the other
component now also contributes to the homology; it is still true that
$\gamma_2$ survives to represent a non-zero element of $H_1(\Gamma_D;
\setQ)$.

When $g=2$, we may simply choose $D$ to consist of a single
non-seperating curve $\gamma$. Then $\Sigma'$ is a connected genus $1$
surface with $r+2$ boundary components. The rational homology of
$\Gamma'$ is thus $\setQ^{r+2}$, and we see that $H_1(\Gamma_D,
\setQ)$ has dimension $r+1$, spanned by the Dehn twists in the
boundary components and $\gamma$.

In the remaining case of $g=1$, first assume $r\geq 1$, and
let $D$ consist of a single curve parallel to a boundary
component. Then obviously $\Gamma_D = \Gamma$, and $H_1(\Gamma_D,
\setQ) = H_1(\Gamma, \setQ) = \setQ^r$. Finally, in the special case
of a closed torus, we refer to the example given in the next section
of a cocycle with values in the module $\mathcal{O}(\mathcal{M})$.

\chapter{Algebraic coefficients}
\label{cha:algebr-coeff}

Although the results in
the present paper indicate that $H^1(\Gamma,
{\mathcal{O}}(\mathcal{M})^*)$ is not trivial, this does not
necessarily imply that the same holds true for the cohomology
$H^1(\Gamma, \mathcal{O}(\mathcal{M}))$ with algebraic functions as
coefficients.

In the simple case of a closed torus, there is an example of a cocycle
with values in $\mathcal{O}(\mathcal{M})$ which cannot be a
coboundary. Namely, consider the well-known presentation of
$\Gamma_{1,0} \cong \SL_2(\setZ)$ as $\langle \tau_\alpha, \tau_\beta
\mid \tau_\alpha\tau_\beta\tau_\alpha =
\tau_\beta\tau_\alpha\tau_\beta, (\tau_\alpha\tau_\beta)^6 =
1\rangle$. The elements $\tau_\alpha$ and $\tau_\beta$ may be realized
as Dehn twists in curves $\alpha$, $\beta$ intersecting transversely
in a single point. We now define a cocycle $u$ on the generators by
$u(\tau_\alpha) = \alpha-\beta$ and $u(\tau_\beta) = \beta - \alpha$,
where on the right hand sides we consider $\alpha$ and $\beta$ as
$1$-component BFK-diagrams. It is easy to check that this in fact
defines a cocycle, since
\begin{align*}
  u(\tau_\alpha\tau_\beta) &= (\alpha-\beta) +
  \tau_\alpha(\beta-\alpha) = -\beta + \tau_\alpha\beta\\
  u(\tau_\beta\tau_\alpha) &= (\beta-\alpha) +
  \tau_\beta(\alpha-\beta) = -\alpha + \tau_\beta \alpha
\end{align*}
so
\begin{align*}
  u(\tau_\beta\tau_\alpha\tau_\beta) &= (\beta-\alpha) +
  \tau_\beta(-\beta + \tau_\alpha\beta) = -\alpha +
  \tau_\beta\tau_\alpha\beta\\
  u(\tau_\alpha\tau_\beta\tau_\alpha) &= (\alpha-\beta)
  +\tau_\alpha(-\alpha + \tau_\beta \alpha) = -\beta +
  \tau_\alpha\tau_\beta\alpha
\end{align*}
But by Lemma~\ref{lem:5}, we have $\tau_\alpha\tau_\beta\alpha =
\beta$ and $\tau_\beta\tau_\alpha\beta = \alpha$, so both right hand
sides are $0$, and $u$ satisfies the first relation. Now it is trivial
to see that it also satisfies the second, because by the first
relation it may also be written $(\tau_\alpha\tau_\beta)^6 =
(\tau_\alpha\tau_\beta\tau_\alpha)^4 = 1$, so we have
\begin{align*}
  u( \psi^4) =  u(\psi) + \psi u(\psi) + \psi^2 u(\psi) = \psi^3
  u(\psi) = 0,
\end{align*}
where $\psi = \tau_\alpha\tau_\beta\tau_\alpha$.

It is clear that $u$ is not a coboundary, because for every linear
combination $f$ of BFK-diagrams, the coefficient of $\alpha$ in
$(\delta f)(\tau_\alpha) = f - \tau_\alpha f$ is necessarily $0$. This
proves that $H^1(\Gamma_{1,0}, \mathcal{O}(\mathcal{M})) \not= 0$. It
is interesting to see if this example can be generalized, for example
using the simple presentation of $\Gamma_{g, r}$ given
in~\cite{MR1851559}.

\appendix

\chapter{Appendix}
\label{sec:appendix}

\section{The moduli space}
\label{sec:moduli-space}

Let $P_i$, $i\in I$, be a collection of pair-wise non-isomorphic
principal $G$-bundles over $\Sigma$, such that any principal
$G$-bundle is isomorphic to some (clearly unique) $P_i$. We let
$\mathcal{A}^F_{P_i}\subset \mathcal{A}_{P_i}$ denote the space of
flat connections. The gauge group $\mathcal{G}_{P_i} = \Aut(P_i)$ acts
on this space, and we let $\mathcal{M}_{P_i} =
\mathcal{A}^F_{P_i}/\mathcal{G}_{P_i}$. We then define the moduli
space of flat $G$-connections over $\Sigma$ to be $\mathcal{M} =
\bigsqcup_{i\in I} \mathcal{M}_{P_i}$.

Choosing a basepoint $x\in \Sigma$, the \emph{representation variety}
is the space
\begin{align}
  \label{eq:4}
  \mathcal{R} = \Hom(\pi_1(\Sigma, x), G)/G,
\end{align}
where the action of $G$ is by post-conjugation. It is well-known that
there is a bijection $R\colon\mathcal{M}\to \mathcal{R}$ given as
follows: For each $i\in I$, choose some $p_i$ in the fibre of $P_i$
over $x$. For a gauge equivalence class $[A]$ of flat connections in
$P_i$ and a homotopy class $[\gamma]$ of loops based at $x$,
$R([A])([\gamma])\in G$ is the holonomy along $\gamma$ with respect to
$A$. This defines a homomorphism $R([A])\colon \pi_1(\Sigma, x)\to G$,
and the dependence on the choice of points $p_i$ vanish when we pass
to the quotient $\Hom(\pi_1(\Sigma, x), G)/G$.

If we choose a finite presentation $\langle a_1, \dots, a_n \mid r_1,
\dots, r_m\rangle$ of $\pi_1(\Sigma, x)$, we may identify
$\Hom(\pi_1(\Sigma, x), G)$ with a certain closed subset $H$ of
$G^{\times n}$, namely the $n$-tuples $(A_1, \dots, A_n)$ such that
$r_i(A_1, \dots, A_n) = 1$ for all $i$. If $G$ is an algebraic group,
we may consider the ring $\mathcal{O}$ of algebraic functions on
$G^{\times n}$, and inside this the ideal $I(H)$ of functions
vanishing on $H$. The ring of algebraic functions on $H$ is the
quotient ring $\mathcal{O}/I(H)$. Since the representation variety is
identified with $H/G$, we define the algebraic functions on
the representation variety to be the set $\mathcal{O}^G/I(H)$; the
space of algebraic functions on $G^n$ which are invariant under
conjugation, modulo the ideal vanishing on $H$. This is, in fact,
independent of the chosen presentation of $\pi_1$, whence there is a
well-defined notion of algebraic functions on the moduli space.

\section{The mapping class group}
\label{sec:mapping-class-group}

The mapping class group $\Gamma$ of an oriented surface $\Sigma$ may
be defined as the group $\Diff(\Sigma, \boundary)$ of
orientation-preserving diffeomorphisms of $\Sigma$ fixing the boundary
pointwise, modulo isotopies fixing the boundary. Since any
homeomorphism is isotopic to a diffeomorphism, and two diffeomorphism
are isotopic through diffeomorphisms if and only if they are isotopic
through homeomorphisms, we may also unambigously speak of the isotopy
class of a homeomorphism of $\Sigma$.

\begin{figure}[htb]
  \centering
  \includegraphics{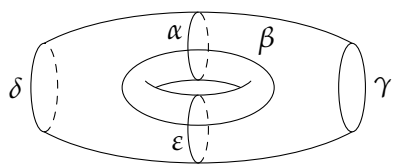}
  \caption{The two-holed torus relation}
  \label{fig:twoholedtorus}
\end{figure}

It is well-known that $\Gamma$ is generated by the isotopy classes of
a certain ``twist diffeomorphisms'' known as Dehn twists, which may be
defined as follows: Let $A$ be the annulus given in polar coordinates
$(r, \theta)$ by $1\leq r\leq 2$, and choose some smooth, increasing
function $\phi\colon [1,2]\to [0,2\pi]$. The standard twist
diffeomorphism of $A$ is given in polar coordinates by $(r,
\theta)\mapsto (r, \theta+\phi(r))$. This diffeomorphism fixes
$\boundary A$ point-wise, and its isotopy class is independent of the
choice of $\phi$. Now, if $\alpha$ is some simple closed curve on
$\Sigma$, we may choose an orientation-preserving embedding of $A$ in
$\Sigma$ such that the inner boundary component coincides with
$\alpha$. Then a twist along $\alpha$ is obtained by extending the
standard twist of $A$ by the identity to the rest of $\Sigma$. (A
priori, this may only define a homeomorphism, but one could for
instance require $\phi$ to be constant near $1$ and $2$). The isotopy
class of this twist is independent of the embedding, and depends only
on the isotopy class of $\alpha$. It is known as the Dehn twist
$\tau_\alpha$ along $\alpha$. We list a few easy facts about Dehn
twists, the first of which is geometrically obvious. For proof we
refer to section~4 of~\cite{MR1886678}.
\begin{lemma}
  \label{lem:3}
  Dehn twists on non-intersecting curves commute. \qed
\end{lemma}
\begin{lemma}
  \label{lem:4}
  If $h$ is a diffeomorphism of $\Sigma$, the
  relation $h\tau_\gamma h\inv = \tau_{h\gamma}$ holds in $\Gamma$.
\end{lemma}
\begin{lemma}
  \label{lem:5}
  If $\alpha$ and $\beta$ are (isotopy classes of) simple closed
  curves intersecting transversely in a single point, then
  $\tau_\alpha\tau_\beta\alpha = \beta$, and the Dehn twists are
  \emph{braided}, ie. satisfy $\tau_\alpha \tau_\beta \tau_\alpha =
  \tau_\beta \tau_\alpha \tau_\beta$.
\end{lemma}
\begin{lemma}
  \label{lem:8}
  Consider the surface $\Sigma_{0,4}$, ie. a sphere with four holes.
  Let $\gamma_i$ denote the $i$'th boundary component, $0\leq i\leq
  3$, and $\gamma_{ij}$ a loop enclosing the $i$'th and $j$'th
  boundary components, $1\leq i<j\leq 3$. Let $\tau_i =
  \tau_{\gamma_i}$ and $\tau_{ij} = \tau_{\gamma_{ij}}$. Then
    \begin{align}
      \label{eq:62}
      \tau_0\tau_1\tau_2\tau_3 = \tau_{12}\tau_{13}\tau_{23}.
    \end{align}
\end{lemma}
This is known as the \emph{lantern relation}.
\begin{proof}
  We may also regard $\Sigma_{0,4}$ as a disc with three open discs
  removed as in Figure~\ref{fig:lantern:0}. Connect $\gamma_0$ with
  $\gamma_i$ by a small arc $I_i$, $i=1,2,3$, such that the three arcs
  are disjoint. Then $\Sigma_{0,4}$ cut along these arcs is a disc
  (with corners), and since the mapping class group of a disc is
  trivial, a diffeomorphism of~$\Sigma_{0,4}$ fixed on $\boundary
  \Sigma_{0,4}$ is determined (up to an isotopy fixed on $\boundary
  \Sigma_{0,4}$) by its action on the arcs $I_1, I_2, I_3$.  Thus one
  need only calculate the effect of both sides of~\eqref{eq:62}
  on~$I_i$ and see that the results agree up to isotopy fixed on
  $\boundary \Sigma_{0,4}$. For simplicity, we only draw the pictures
  relevant for $I_1$; the reader can easily draw the corresponding
  pictures for the other two arcs.
  \begin{figure}[htb]
    \centering
    \hfill
    \subfloat[\label{fig:lantern:0} A sphere with four
    holes.]{\includegraphics{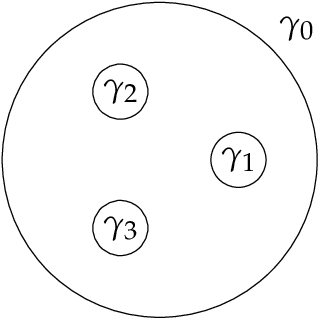}}
    \hfill
    \subfloat[\label{fig:lantern:1} Loops enclosing boundary
    components.]{\includegraphics{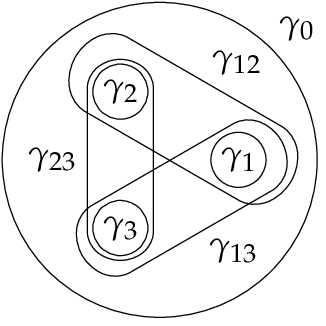}}
    \hfill
    \subfloat[\label{fig:lantern:2} Arcs connecting boundary
    components.]{\includegraphics{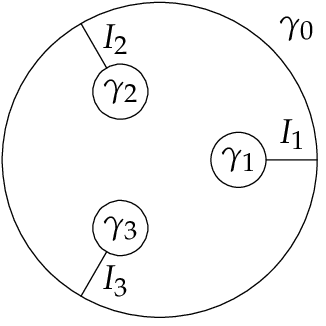}}
    \hfill\strut
    \caption{The lantern relation.}
    \label{fig:lantern}
  \end{figure}

  Now if we choose $I_1$ to be the horizontal line segment of
  Figure~\ref{fig:lantern:2}, we see that $\tau_3$, $\tau_2$, and
  $\tau_{23}$ act trivially on $I_1$. Thus we need only to show that
  $\tau_0\tau_1$ and $\tau_{12}\tau_{13}$ has the same effect on
  $I_1$. This is clear from the two rows of pictures in
  Figure~\ref{fig:lanternproof} below.
\end{proof}
\begin{figure}[htb]
  \centering
  \settowidth\unitlength{$\scriptstyle \tau_{13}$}
  \newcommand\ma[1]{\xrightarrow{\mathmakebox[\unitlength][c]{#1}}}
  $\vchb{\includegraphics{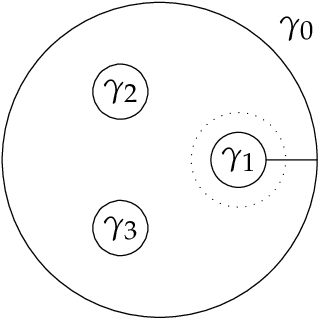}}\hfill\ma{\tau_1}\hfill
   \vchb{\includegraphics{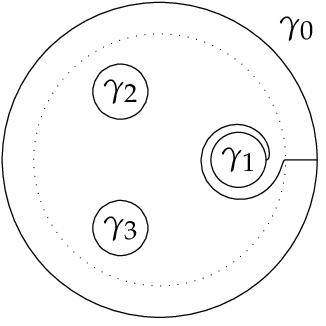}}\hfill\ma{\tau_0}\hfill
   \vchb{\includegraphics{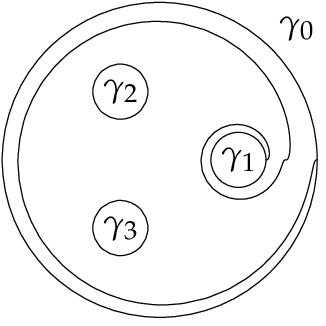}}$%
  \medskip

  $\vchb{\includegraphics{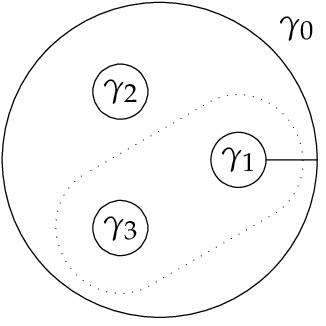}}\hfill\ma{\tau_{13}}\hfill
   \vchb{\includegraphics{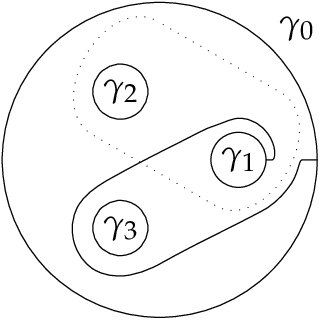}}\hfill\ma{\tau_{12}}\hfill
   \vchb{\includegraphics{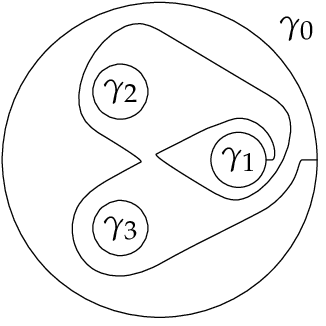}}$
  \caption{The effect of $\tau_0\tau_1$ and $\tau_{12}\tau_{13}$ on
    $I_1$ agree up to isotopy.}
  \label{fig:lanternproof}
\end{figure}
\begin{lemma}
  \label{lem:6}
  Let $\Sigma_{1,2}$ be a two-holed torus, and let $\alpha, \epsilon$ be
  non-intersecting simple closed curves which both intersect a simple
  closed curve $\beta$ in a single point. Further, let $\delta,
  \gamma$ denote simple closed curves parallel to the boundary
  components (see Figure~\ref{fig:twoholedtorus}). Then
  \begin{align}
    \label{eq:5}
    (\tau_\alpha\tau_\epsilon\tau_\beta)^4 = \tau_\delta \tau_\gamma.
  \end{align}
\end{lemma}
This lemma can be proved in a completely similar fashion as above by
choosing a set of proper arcs such that cutting along these gives a
disc (with corners), and proving that both sides of \eqref{eq:5} has
the same effect on this set of arcs.

These simple results have a few interesting consequences with respect
to the homology of the mapping class group
\begin{proposition}
  \label{prop:3}
  If $g\geq 2$, the Dehn twist on a boundary component of $\Sigma_{g,r}$
  can be written in terms of Dehn twists on non-separating curves.
\end{proposition}
\begin{proof}
  The assumption on the genus implies that we may find an embedding of
  $\Sigma_{0,4}\to\Sigma_{g, r}$ such that $\gamma_0$ is mapped to the
  boundary component in question and the remaining six curves involved
  in the lantern relation are mapped to non-separating curves (think
  of $\Sigma_{g,r}$ as being obtained by gluing three boundary
  components of $\Sigma_{g-2, r+2}$ to $\gamma_1$, $\gamma_2$ and
  $\gamma_3$, respectively). Then the relation $\tau_0 =
  \tau_{12}\tau_{13}\tau_{23}\tau_3\inv\tau_2\inv\tau_1\inv$ also
  holds in $\Gamma_{g,r}$.
\end{proof}
\begin{proposition}
  \label{prop:4}
  When $g\geq 2$, $\Gamma_{g,r}$ is generated by Dehn twists on
  \emph{non-separating} curves.
\end{proposition}
\begin{proof}
  We already know that the mapping class group is generated by Dehn
  twists, so it suffices to show that a Dehn twist on a seperating
  curve $\gamma$ can be written in terms of twists on non-separating
  curves in $\Sigma$. If we assume $g\geq 3$, cut $\Sigma$
  along $\gamma$ and apply Proposition~\ref{prop:3} to the component
  which has genus $\geq 2$.

  Now, if $g = 2$ and $\gamma$ seperates $\Sigma$ into a genus $0$ and
  a genus $2$ component, we may still apply Proposition~\ref{prop:3}
  to see that $\tau_\gamma$ is a product of twists on non-seperating
  curves. Hence assume $\gamma$ cuts $\Sigma$ into two genus $1$
  components. If one of these components has a boundary
  component other than the one arising from the cut along $\gamma$, we may
  apply Lemma~\ref{lem:6} to see that $\tau_\gamma$ can be written in
  terms of twists on non-seperating curves along with the twist on
  the additional boundary component. But the latter may be written in
  terms of twists on non-seperating curves in the original surface. If
  there are no other boundary components, we may still apply
  Lemma~\ref{lem:6}, since we can simply cut out a disc bounded by a
  trivial loop $\delta$; then the relation~\eqref{eq:5} degenerates to
  $(\tau_\alpha\tau_\epsilon\tau_\beta)^4 = \tau_\gamma$.
\end{proof}
\begin{corollary}
  \label{cl:1}
  When $g\geq 3$, $H_1(\Gamma, \setZ) = 0$, and when $g = 2$,
  $H_1(\Gamma, \setZ) = \setZ/10\setZ$. In both cases, $H_1(\Gamma,
  \setQ) = 0$.
\end{corollary}
\begin{proof}
  Assume $g\geq 2$. The first homology group with integral
  coefficients is the same as the abelianization of the group. By
  Proposition~\ref{prop:4}, $\Gamma$ is generated by Dehn twists on
  non-seperating curves. Since any two non-seperating curves are
  related by a diffeomorphism, such Dehn twists are always conjugate
  by Lemma~\ref{lem:4}. Hence $H_1(\Gamma, \setZ)$ is cyclic,
  generated by any such Dehn twist~$\tau$. If $g\geq 3$, we may find
  an embedding of $\Sigma_{0,4}$ into $\Sigma$ such that all seven
  curves occuring in the lantern relation~\eqref{eq:62} are
  non-seperating, so the relation $4\tau = 3\tau$ holds in the
  abelianization, and $H_1(\Gamma, \setZ) = 0$ when $g\geq 3$.

  If $g = 2$, there is an embedding of the two-holed torus in $\Sigma$
  such that all five curves occuring in Lemma~\ref{lem:6} are
  non-seperating. This implies that the relation $12\tau = 2\tau$
  holds in $H_1(\Gamma, \setZ)$; or in other words that $H_1(\Gamma,
  \setZ)$ is cyclic of order dividing $10$. To see that the order is
  in fact $10$, we observe that there is an epimorphism $\Gamma_{2, r}
  \to \Gamma_2$ inducing an epimorphism of abelianizations
  $H_1(\Gamma_{2, r}, \setZ)\to H_1(\Gamma_2, \setZ)$, so it suffices
  to prove the statement for a closed surface of genus $2$. The
  easiest way to do this is to use one of the known presentations of
  the mapping class group; for instance the very symmetric
  presentation given by S.~Gervais in~\cite{MR1851559}.
\end{proof}
The first proof that $H_1(\Gamma, \setZ) \cong \setZ/10\setZ$ is due
to Mumford~\cite{MR0219543}, and does not rely on a presentation of
the mapping class group. We also remark that for $g=1$ and $r\geq 1$,
the result in Proposition~\ref{prop:2}(\ref{genus-eq-1}) can also be
obtained from the presentation given in~\cite{MR1851559}.

\bibliographystyle{plainurl}
\bibliography{../phd}

\end{document}